\documentclass{ieeeaccess}
\usepackage{siunitx}
\usepackage{amsmath,amssymb,amsfonts}
\usepackage{algorithmic}
\usepackage{graphicx}
\usepackage{textcomp}
\def\BibTeX{{\rm B\kern-.05em{\sc i\kern-.025em b}\kern-.08em
    T\kern-.1667em\lower.7ex\hbox{E}\kern-.125emX}}
\begin{document}
\history{Date of publication xxxx 00, 0000, date of current version xxxx 00, 0000.}
\doi{xx.xxxx/ACCESS.2023.DOI}

\title{Power Reduction in FM Networks by Mixed-Integer Programming. A Case Study}
\author{\uppercase{Pasquale Avella}\authorrefmark{1},
\uppercase{Paolo Nobili\authorrefmark{2}, and Antonio Sassano}.\authorrefmark{3}}
\address[1]{Dipartimento di Ingegneria, Università del Sannio, Piazza Roma 21, 82100 Benevento Italy (e-mail: avella@unisannio.it)}
\address[2]{DAFNE, Università della Tuscia, Via San Camillo de Lellis 01100 Viterbo Italy(e-mail: nobili@unitus.it)}
\address[3]{Dipartimento di Ingegneria Informatica, Automatica e Gestionale - DIAG, Sapienza Università di Roma, via Ariosto 25, 00185 Roma Italy (email: sassano@diag.uniroma1.it}

\corresp{Corresponding author: Pasquale Avella (e-mail: avella@unisannio.it).}

\begin{abstract}
\noindent The climate change emergency calls for a reduction in energy consumption in all human activities and production processes. The radio broadcasting industry is no exception. However, reducing energy requirements by uniformly cutting the radiated power at every transmitter can potentially impair the quality of service. A careful evaluation and optimization study are in order. In this paper, by analyzing the Italian frequency modulation analog broadcasting service, we show that it is indeed possible to significantly reduce the energy consumption of the broadcasters without sacrificing the quality of the service, rather, even getting improvements.
\end{abstract}

\begin{keywords}
Broadcast networks, FM, Power reduction, Mixed-Integer Programming
\end{keywords}

\maketitle

\sisetup{group-separator={,}}

\section{Introduction}\label {sec:Introduction}

\noindent The frequency-modulated analog radio (FM) service has been, for over 60 years, a faithful companion on long car journeys and, for many, the only alternative to loneliness in their homes. Now, the end of its natural life cycle is approaching: many European countries (Switzerland, UK, Denmark, Sweden, Germany) have devised and implemented plans for its definitive shutdown; Norway already dismissed all national networks in 2017. However, FM analog broadcasting still resists digital radio and Internet streaming replacements and we can safely predict that it will be active in several countries, including Italy, for quite a few years to come. 

\medskip

\noindent
In fact, all shutdown projects have met with fierce resistance from users and (mostly local) broadcasters. Users do not want to abandon the FM. However, the reasons for the switch-off were and are very valid. Digital radio and, above all, Internet streaming are much more efficient and widespread vectors for radio programs. Thus an Italian user (for example) can listen to his favorite local radio in Denmark without any problems. However, not everyone has a digital receiver or a smartphone or, rather, not everyone is used to turn on a smartphone to listen to radio broadcasting, and not all cars are equipped with DAB receivers (only since 2018 has the presence of such receivers become mandatory for new cars in Italy). Therefore, technological evolution alone has not been able to change user's habits.

\medskip

\noindent To this technological evolution, more recently the need (by administrations) has been added to avoid wasting energy and hope (by broadcasters) to stay in business in the face of the rising costs of electricity bills. To this end, many are asking to consider the FM broadcasting industry an energy-intensive sector and request state aid to continue broadcasting. However, is it an energy-intensive sector? In 2020 the British Broadcasting Company (BBC) \cite{FletcherBBC} published a study entitled "The energy footprint of BBC radio services: now and in the future" which evaluated the energy consumption associated with the broadcasting and receiving of radio programs. The results were astonishing: the energy requirements of the analog component (essentially only the FM radio) computed considering both the receivers and the broadcasting infrastructure, unquestionably emerged as the most demanding, as shown in Table \ref{tab:BBC}.

\begin{table}[hhh] 
\caption{Top five components with the highest energy consumption in the BBC radio system for 2018}
\label{tab:BBC}
\begin{tabular}{|l|c|c|}
\hline
Radio System Component      & Annual Energy (GWh) & \% of total \\
\hline
Analogue Radio Set          & 82                  & 25.2                     \\
DAB Radio Set               & 55                  & 17.0                     \\
TV Set                      & 33                  & 10.1                     \\
FM Broadcast Infrastructure & 26                  & 8.1                      \\
AM Broadcast Infrastructure & 17                  & 5.2                      \\
\hline
\end{tabular}
\end{table}

\bigskip\noindent
The BBC estimates the total energy consumed by all forms of radio and TV broadcasting in 2018 at $325$ GWh, corresponding to a $37.1$ MW stove lit for $24$ hours a day and $365$ days a year. Moreover, $100$ GWh ($30.8\%$) was used by the FM service alone, of which approximately $26$ GWh ($8.1\%$) by the broadcast infrastructure. Analog (FM) is the most energy-intensive method for broadcasting radio programs. Hence the urge to “switch-off” analog broadcasting which, unfortunately, would have, at the moment, unsustainable social costs.

\medskip\noindent
Usually, the problem of reducing the power of FM networks is addressed by designing more efficient directional antennas \cite{Mappatao1}. Instead, an optimization-based approach would allow a reduction in energy consumption without decreasing the quality of the service, being aware that at present there is a strong waste of resources. Obviously, this solution does not have the same effect for all countries but it is certainly more effective in countries where the "waste" is more significant, and Italy is one of them.

\medskip

\noindent
In the UK the BBC broadcasts $40$ programs while there are approximately 450  other local and national broadcasters; there are approximately 1,500 transmitters; In France, approximately 2,200 transmitters broadcast $945$ programs. In Italy, according to the FM Frequency Register developed by the State Agency for Telecommunications (AGCOM), there are approximately 16,381 transmitters broadcasting programs of 1,133 networks \cite{FM Frequency Register}. These numbers alone do not say much about the energy consumption. A few very powerful transmitters could consume the same energy as hundreds of small transmitters. However, Italy has almost eight times more transmitters than France, compared with a slightly higher number of broadcasters (1,133 vs.\ 945). The number of frequencies used by Italian networks was at least twice that used by French networks. This situation tends to produce a significant "national self-interference" (i.e. Italian broadcasters disturbing other Italian broadcasters) which, in turn, is a push to increase the transmission power to counteract the interference. 
\medskip

\noindent In short, there are several signs supporting the hypothesis that the average electromagnetic radiated power (and hence the power consumed) by transmission plants in Italy may be particularly high. In recent years, specialized online magazines have dealt with this problem seriously, underlining how (for example) the power of the FM transmitters in the city of Milan has had to increase by $10$ times to counter the interference generated by other (Milanese) broadcasters and how, due to the age of the transmission plants, the efficiency of the systems (the ratio between the electromagnetic power radiated at the antenna and the power consumed at the meter) is far from the ideal $70\%$.

\medskip\noindent Finally, the FM broadcasters themselves declare that they manage a network of transmitters that consume a large amount of energy and that they are forced to decrease the power of their networks during the night to cut energy costs. In short, everything would suggest that the Italian FM "industry" is certainly energy-intensive but also particularly inefficient in ensuring service to users.

\medskip\noindent To quantify the energy consumption of the 16,381 Italian FM transmitters, we used the FM frequency register developed by AGCOM. The report only provides the amount of radiated power from each of the Italian plants, but we can infer a good approximation of the energy consumption by guessing the (mean) efficiency of the transmitters. Given the old age of most Italian plants, a good estimate for such efficiency would be approximately $50\%$, and we used such a hypothesis in our computations. The result of the calculation is that the energy consumption of Italian FM plants is approximately $253$ GWh per year (in contrast to $26$ GWh reported by the BBC for the United Kingdom - see above) and with an average power of approximately 29 MW. In short, the consumption of Italian FM networks alone is close to the total consumption of British radio networks and is approximately 10 times higher than that of British FM networks alone.

\section{Literature}\label{sec:Literature}

The scientific literature concerning planning and re-planning of Wireless Networks with the goal of energy sustainability has flourished in recent years. \textit{5G Networks}, \textit{sensor networks} and the \textit{Internet of Things (IOT) environment} ask for accurate service and low energy consumption. Consequently the problem of maximizing the energy efficiency of the system by adjusting the minimum signal-to-interference plus noise ratio to guarantee the required service has been studied by many authors see \cite{Mavrommatis},  \cite{gonzalez2011base}, and \cite{ParkLim}. With the exception of \cite{Mappatao1}, the FM Networks have received much less attention and our paper, to the best of our knowledge, is the first to pose in an optimization format the problem of re-planning an FM network to reduce energy consumption, preserve and improve the current service and reducing the interference produced in the bordering Administrations.

\medskip

\section{Definitions and a practical example}\label{sec:Definitions}

\noindent The FM infrastructure is a set of transmitters $T$ broadcasting radio programs to a set $R$ of geographically distributed locations, represented by \emph{receiving points} on a map. Let $A$ denote the set of the radio broadcasting networks. We assume that transmitter $t \in T$ carries the programs of a single network $a \in A$, being $T(a)$ the set of transmitters of the network $a \in A$. Transmitter $t$ operates on a single frequency and emits a radio signal with power $p_t \in [0, P_{max}]$. The power $p_{rt}$ received at a specific geographic location $r$ from transmitter $t$ is proportional to the emitted power $p_t$ by a factor $a_{rt} \in [0, 1]$, that is $p_{rt}  = a_{rt} \cdot p_t$. The factor ${a}_{rt}$ is called the \emph{fading coefficient} and summarizes the reduction in power that a signal experiences while propagating from $t$ to $r$. The value of the fading coefficient depends on many factors (e.g., the distance between the communicating devices, presence of obstacles, antenna patterns, and tropospheric effects) and is commonly computed using a suitable propagation model. 

\medskip\noindent
Let $T(r)$ denote the set of transmitters received at point $r$. To be acceptably listened to in $r$, the received power must be greater than a minimum value $p_{min}$, called \textit{background noise}, and overcome the \textit{interference} of all the other transmitters received in $r$. Note that a transmitter $t$ can also interfere in $r$ if $p_{rt} < p_{min}$. Usually, the interfering power of a transmitter in $r$ is computed by considering a different fading factor ${\bar a}_{rt}$ and is multiplied by a constant called the protection ratio $PR$, namely $\bar p_{rt}  = \bar a_{rt} \cdot PR \cdot p_t$. Hence, for each transmitter $t$ and receiver $r$ we have two different receiver powers, namely, the \textit{useful power} $p_{rt}$ and the \textit{interfering power} $\bar p_{rt}$ with, usually, $p_{rt} \le \bar p_{rt}$. The effect of interfering signals on a useful signal is a complex physical phenomenon that determines the reception quality. 

\medskip\noindent
The interference is stronger when useful and interfering signals are transmitted at the same frequency (\textit{co-channel interference}) but it is non-negligible also if the signals are modulated at different frequencies. In this paper, we will consider only co-channel interference and we will denote by $I(rt)$ the set of transmitters, different from $t$, received in $r$ at the same frequency used by transmitter $t$. Following a standard technical practice \cite{Rappaport}  we assume that a useful signal from $t$ can be received with acceptable quality at the receiving point $r$ if the ratio between the useful signal and the sum of the interfering signals and background noise, called \emph{signal-to-interference-plus-noise ratio} (\textit{SINR}) is above a threshold $\theta$:
\begin{equation}
  \frac{p_{rt}}{\sum_{j \in I(rt)}  \bar p_{rj} +  p_{min}} \ge \theta
\end{equation}

\smallskip\noindent It is important to note that the useful signal is emitted by a unique transmitter, and all the other co-channel signals interfere. This is due to the fact that Frequency Modulation is an \textit{analog} technology and, contrary to the more advanced digital technologies like DAB, there is no positive composition of co-channel signals carrying the same content.

\medskip\noindent 
A comment about our choice of studying the case in which the interference is produced only by transmitters operating at the same frequency (\textit{Co-Channel Hypothesis}) is mandatory here. This choice is strongly motivated by the main goal of our paper, namely showing that the reduction of the operating powers of the FM transmitters preserves the current service while reducing the energy consumption and the Italian interference in the bordering Administrations. The question is: How do we define the "current service"? A rigorous definition should take into account the interference generated by the adjacent channels but, in this case, we would reduce the current service areas of the Italian broadcasters and hence we would favor a more drastic reduction in operating powers. The Italian broadcasters would certainly object that their service is "acceptable" even in the presence of adjacent channel interference and this claim is strongly supported by the empirical evidence that the Italian FM operating frequencies are seldom separated by the $200$ KHz suggested by ITU \cite{ITU_FM}. This has been the main motivation for adopting the Co-Channel Hypothesis: \textit{the maximization of the current service areas of Italian networks to protect}. 

\medskip\noindent Another important point underlines the importance of the choice of the Italian networks in our case study. In general, many transmitters of the same network, possibly modulated at different frequencies, are received at the same receiving point $r$. This is usually a clear sign of a network that has been poorly designed or that has grown in a "tragedy of the commons" environment. In a healthy environment in which frequencies are properly used and accurately planned, the average number of transmitters received from each network at each receiving point should be close to one. This is the case for most FM networks in the world but not in Italy, making the choice of our case study particularly interesting. It follows that in our study we can have more than one transmitter $t$ carrying the same network and with a SINR above the threshold. 

\medskip\noindent 
Let $P_r$ be the population of the receiver point $r$. Let $y_t$ be the percentage reduction in the power of transmitter $t \in T$. For each receiving point $r \in R$ and each network $a \in A$, let $s_{ra}$ be a binary variable that is $0$ if the best server of network $a$ in $r$ has a SINR above the threshold (covered) and $1$ otherwise. Set $T$ contains only the servers of some receiving point $r\in R$. Let $t_{ra} \in T(a)$ be the best server of the network $A$ for the receiving point $r$ and let $p_{ra}$ be the useful power received in $r$ from $t_{ra}$. Let $I_{ra}$ be the set of co-channel interfering transmitters for $t_{ra}$ and $\theta$ the threshold above which a receiver point is regarded as covered. Finally, let $Z$ be the set of pairs $(r,a)$ with the property that the server of the Italian network $a$ in an Italian receiving point $r$ satisfies the service constraint (SINR ratio) in the current configuration. That is $Z$ represents the current FM service to be preserved. 

\medskip\noindent
The task of planning (and re-planning) FM networks is usually performed using powerful and effective \textit{simulation tools}. In a simulation environment, some of the restrictions we are forced to introduce into our optimization model can be relaxed, and the computation can be made more accurate. One of the features of the simulation models is that the quality of service and the level of interference can be presented through informative service maps. In this paper we use a simulation tool to verify the quality of the solutions produced by our optimization models. Before defining the models, we will use our simulation tool to show what happens at a specific receiving point, at a specific frequency when the power of the main interfering transmitter is reduced. 

\medskip\noindent Figures \ref{fig1} and \ref{fig2} focus on the improvement obtained at a real receiving point situated in the vicinity of the Slovenian city of Capodistria. The tables under the pictures report the relevant data regarding the seven transmitters whose signals are received at this point, namely the identifier, \textit{useful power} (in dB), and  \textit{interfering power} (in dB). Of the seven transmitters, the first one is Slovenian (and is the server for this point of a Slovenian network) and the remaining ones are Italian (and hence are interfering). Column \textit{SumInt} reports the cumulative interference provided by the background noise $p_{min}$ and by the interfering transmitters listed in the preceding columns (note that values expressed in dB are not additive).

\Figure[hhh](topskip=0pt, botskip=0pt, midskip=0pt)[width=0.9\linewidth]{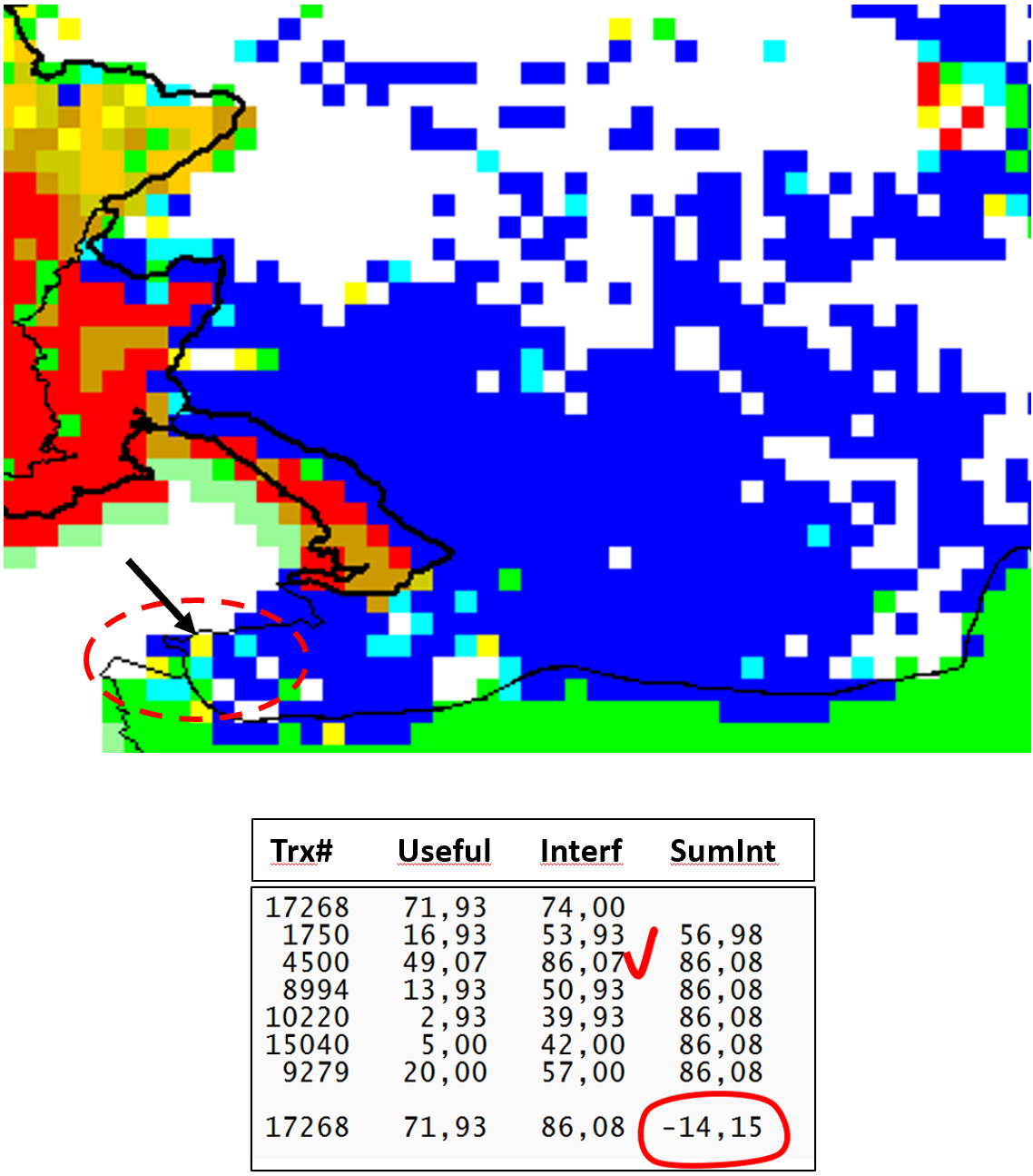} {Capodistria: Current Scenario.\label{fig1}}

\Figure[hhh](topskip=0pt, botskip=0pt, midskip=0pt)[width=0.9\linewidth]{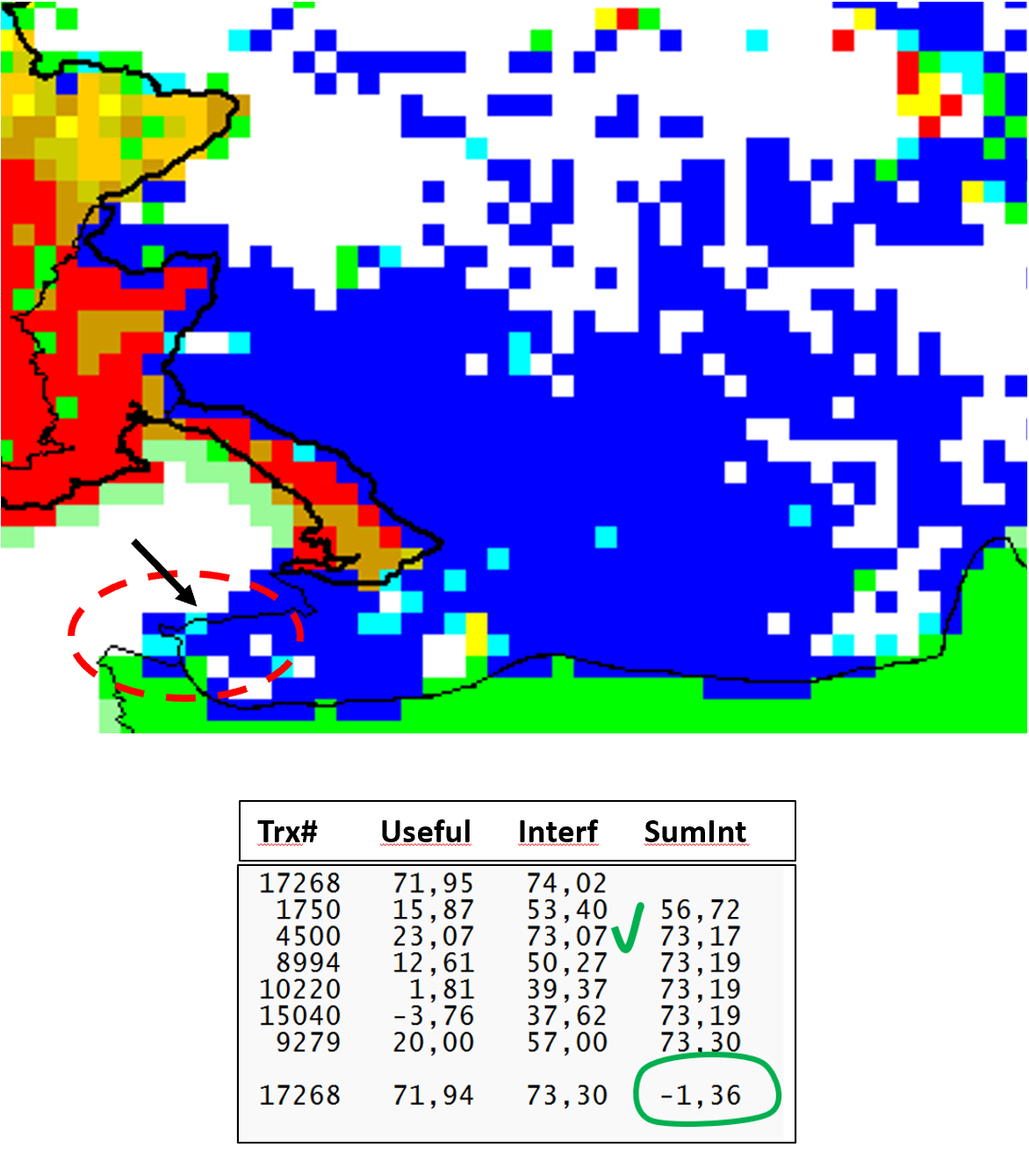} {Capodistria; Post Power Reduction.\label{fig2}}

\medskip\noindent To read the two service maps better we now explain our color code. We start by noticing that the "pixels" in the map are the squares clearly visible at the resolution we have chosen for our example. We call \textit{Affected Administration} the Administration whose Quality of Service (QoS) we want to assess and \textit{Interfering Administration} the Administration whose networks disturb the service of the Affected Administration. Each  "pixel" of the map belonging to the Affected Administration has a color indicating the  QoS of the corresponding receiving point. The possible colors are \textit{(blue, light blue, green, yellow, and red)}. The first four colors are associated with four QoS levels, \textit{(Q4, Q3, Q2, Q1)}. Each grade on the quality scale is achieved when the SINR at the receiving point is greater than or equal to a given threshold. In our case, the thresholds were  $(0, -6, -12, -15)$ \textit{dB}. Finally a red "pixel" indicates that the SINR at the corresponding receiving point is below a threshold of $-18$ dB and no service is available.

\medskip\noindent For the Interfering Administrations, the colors of the "pixels" \textit{(red, brown, light brown, orange, yellow, green)} have different meanings. They represent the level of cumulative interference produced by the networks belonging to the Interfering Administration and are determined by the value of the electric field at the receiving point. In our case, the color changes from \textit{red} (severest interference) to \textit{yellow} (light interference) when the electric field is greater than or equal to the thresholds $(70, 50, 40, 30, 20)$ dB($\mu V/m$). The green color indicates that the signal was received with an electric field below $20$ dB($\mu V/m$) causing negligible interference. 

\medskip\noindent Let us now return to our example and observe that, in this case, Slovenia is the Affected Administration and Italy the Interfering Administration. The interfering power of the Italian transmitter 4500 before the optimization is 86.07 dB, making its signal the worst interfering transmitter received at this point. The cumulative interference experienced by the point is at $86.08$ dB, overcoming the useful power of the Slovenian transmitter (the server) by $14.15$ dB (see the last row in the table) and rendering its reception very poor (i.e.\ the SINR of the Slovenian transmitter at the point is $-14.15$ dB, corresponding to a QoS of \textit{Q1}). In contrast, after optimization, the interfering power of the interfering transmitter $4500$ is decreased to $73.07$ dB and consequently, the cumulative interference is lowered to $73.30$ dB. Hence, after optimization, the SINR of the Slovenian transmitter increases to $-1.36$ dB, corresponding to a QoS of Q3.

\medskip\noindent
The question now is: \textit{What happens to the service of the Italian transmitter 4500 in Italy?} Reducing its power could make some other Italian transmitter, operating at the same frequency, the best server in some important receiving point in Italy. Hence the network operating transmitter 4500 could lose an important portion of its users. In such a situation, the power of the new best server should inevitably be decreased. We would have the  "butterfly effect" phenomenon, well known to all the practitioners trying to resolve interference problems by means of simulation tools, one transmitter at a time. The answer is very simple, the power of $16,381$ Italian transmitters \textit{must be changed simultaneously}, by solving a global optimization problem that guarantees each broadcaster to preserve its original service. This is exactly what we have done in the next section.

\section{The MILP model}\label{sec:MILP model}

In this section, we introduce the optimization models aimed at solving the \textit{FM Power Reduction problem}. In particular, we begin with the following mixed-integer fractional programming problem:
\begin{equation} \label{FracLP}
\begin{aligned}
\min & \sum_{r \in R}P_r\sum_{a \in A}s_{ra} \\
      & \frac{p_{rt_{ra}}y_{t_{ra}}}{\sum_{j \in I_{ra}}\bar{p}_{rj} + p_{min}} + M s_{ra}\ge \theta & r \in R,a\in A \\
      & s_{ra}\in \{0,1\} \ r \in R, a \in A \\
      & s_{ra} = 0, \ (r,a) \in Z \\
0 \le & y_i \le 1, \ i \in T 
\end{aligned}
\end{equation}

\smallskip\noindent
where 
\begin{itemize}
    \item[-] The objective function aims to minimize the total population not covered by FM networks.
    \item[-] $M$ is a very large value that enforces the SINR inequality associated with the pair $(r,a)$ to satisfy exactly one of the following conditions: \textit{either} $i)$ the receiver $r\in R$ is served by the transmitter $t_{ra}$; \textit{or} $ii)$ $s_{ra}$ is set to $1$ and the inequality is no longer active.  
    \item[-] The constraint $s_{ra} = 0$ for $(r, a) \in Z$ forces server $a\in A$ to satisfy the SINR constraint in $r \in R$ and hence to preserve the original service of the network $a$ in $r$.
\end{itemize}
\medskip\noindent 
The above model can be easily linearized as follows:
\begin{equation}
\label{Model MILP}
\begin{aligned}
\min & \sum_{r \in R}P_r\sum_{a \in A}s_{ra}  \qquad (MILP) \\
      & y_{t_{ra}}-\theta\sum_{j \in I_{ra}}\frac{\bar p_{rj}}{p_{rt_{ra}}} y_j \\
      & \qquad + M s_{ra}\ge \frac{p_{min}}{p_{rt_{ra}}}\theta \quad  r \in R, a\in A \\
      & s_{ra}\in \{0,1\} \ r \in R, a \in A \\
      & s_{ra} = 0, \ (r,a) \in Z \\
0 \le & y_i \le 1, \ i \in T
\end{aligned}
\end{equation}

\smallskip \noindent This choice has a positive algorithmic outcome because each SINR inequality involves only transmitters operating at the same frequency, and the model is decomposable into blocks, one for each frequency in the spectrum. This property makes the model more easily solvable by the current MIP algorithms, which can recognize blocks and solve each of them as a sub-MIP.

\medskip\noindent Usually, the LP-relaxation of a MILP containing a Big-M coefficient does not provide good lower bounds and very rarely provides good quality feasible solutions to the original problem. To contrast this bad behavior of the Big-M constant and produce more efficient lower and upper bounds, the authors of \cite{avella2023} used a different (and provably lower) value for the Big-M. Namely, the sum of all the power received from the interfering transmitters and the background noise or, more formally, $M_{ra}=\theta\frac{p_{min}}{p_{rt_{ra}}}+ \theta\sum_{j \in I_{ra}}\frac{\bar p_{rj}}{p_{rt_{ra}}}$. 

\medskip\noindent It is expected that this lower value of the constant $M_{ra}$ could produce a better relaxation and more accurate upper and lower bounds. One of the surprising results of this paper is that this does not happen and the results produced by the formulation in \cite{avella2023} lead to power reductions and FM service of inferior quality with respect to those obtained in our MILP model with $M = 10^{40}$. We have observed above that this behavior could be motivated by the efficient management of the Big-M constraints ensured by Gurobi (while the lower values of $M_{ra}$ and their high variability could have a masquerading effect for Gurobi). However, what is really surprising and promising for all practitioners interested in replicating our results is that a special LP problem which is not the relaxation of our original MILP produces results, in terms of power reduction and service increase, which are even better than those produced by the MILP model. The LP model is the following:

\begin{equation} \label{Model LP}
\begin{aligned}
\min & \sum_{r \in R}P_r\sum_{a \in A}s_{ra} \qquad (LP)\\
      & y_{t_{ra}}-\theta\sum_{j \in I_{ra}}\frac{\bar p_{rj}}{p_{t_{ra}}} y_j \\
      & \qquad + s_{ra}\ge \frac{p_{min}}{p_{rt_{ra}}}\theta \quad  r \in R, a\in A \\
      & s_{ra} \ge 0 \ r \in R, a \in A \\
      & s_{ra} = 0 \ (r,a) \in Z \\
0 \le & y_i \le 1, \ i \in T\\
\end{aligned}
\end{equation}

\smallskip\noindent
In this model, the variable $s_{ra}$ can achieve any positive value, and the numerical problems originating by $M$ (or, worst, $M_{ra}$) are overcome. The first advantage of the model (\ref{Model LP}) is its \textit{scalability}; namely, it can be easily solved by any commercial LP-solver in a very short time even if the number of transmitters and receiving points becomes very large. The second - surprising - advantage is that the quality of the results produced by the LP model and in particular the values of the variables $y$, guarantee, as we will see in the next section, a greater reduction in the transmitting powers and a more significant growth in the service area of the Italian networks.  

\section{Results}\label{sec:results}

\medskip

\noindent The MILP (\ref{Model MILP}) and LP (\ref{Model LP}) models were embedded via a Python interface into the MIP solver Gurobi 10 \cite{Gurobi} to obtain provably good solutions. The experiments were run on an Intel Xeon W-10885 $2.40$ GHz workstation with $128$ GB RAM.

\noindent The test bed involves all the Italian transmitters listed in the FM Register managed by AGCOM \cite{FM Frequency Register} and all the foreign FM transmitters registered in the \textit{ITU Master Register of Geneva}, whose service areas are interfered by the Italian FM networks. In particular, all the information was updated on \textit{May 10, 2022}. 

\noindent When talking about service we refer to the service of a transmitter of Administration A on a receiving point located in the territory of Administration A (e.g. the service of an Italian network in Italy). As stated above our goal is to reduce the power of the Italian transmitters (the transmitters belonging to \textit{non-Italian} administrations are fixed to their current power) while guaranteeing that the service areas of the Italian networks are preserved (at least one server of the network must guarantee the service at a receiving point RP if RP is served in the current served in the \textit{current scenario}). The current scenario is summarized in Table \ref{table:current}, where \textit{Population} is the sum over all the radio networks of the users that could be potentially served if there was no interference and \textit{Currently served population} denotes the sum of the users currently covered by all the networks.

\begin{table}[hhh]
\caption{Current scenario} \label{table:current}
\setlength{\tabcolsep}{3pt}
\begin{tabular}{|l|r|}
\hline
Total number of transmitters& 21,805 \\
Number of servers& 20,558 \\
Number of Italian servers& 16,381 \\
Number of foreign servers& 4,177 \\
Population - Italy& 1,834,629,759 \\
Population - abroad& 202,050,318 \\
Currently served population - Italy& 1,377,115,927 \\
Currently served population - abroad& 169,894,991 \\
\hline
\end{tabular}
\end{table}

\noindent With several thousands of active transmitters and only about $200$ channels to share in the FM spectrum, some level of co-channel interference is unavoidable. What is more peculiar to the Italian context (and a sign of poor network design) is that several receiving points experience co-channel interference even from transmitters associated with the same network.

\noindent With the MILP model (\ref{Model MILP}), Gurobi provided in 4,200 secs and 758,237 branch-and-bound nodes, a lower bound LB and an upper bound UB for the optimal solution, with a percentage gap (computed as $100*(UB-LB)/LB$) less than 1\%. Before starting enumeration, the lower bound given by LP-relaxation was strengthened by the Gurobi built-in cutting planes. The outcomes for the MILP model (\ref{Model MILP}) are reported in Table \ref{tab:MILP}, where \textit{Number of Italian plants shut down} counts the transmitters with zero potentially serviceable points (the receiving points in which the transmitter is above $p_{min}$), \textit{$\Delta$ Power} is the percentage variation (decrease in this case) of total transmission power and \textit{$\Delta$ Served population} is the variation (increase) in the sum of served users for all the networks, distinct in Italian and non-Italian but belonging to one of the aforementioned bordering Administrations. 

\begin{table}[hhh] 
\caption{Summary of the outcomes of the MILP model} \label{tab:MILP}
\setlength{\tabcolsep}{3pt}
\begin{tabular}{|l|r|}
\hline
Number of Italian plants shut down& 1,473 \\
$\Delta$ Power (\%)& -53.23 \\
$\Delta$ Served population - Italy& +133,824,096 \\
$\Delta$ Served population - abroad& +16,531,361 \\
\hline
\end{tabular}
\end{table}

\noindent Gurobi solved to optimality the LP model $(\ref{Model LP})$ in only $21$ seconds. The outcomes of the LP model are presented in Table \ref{tab:LP}. 

\begin{table}[hhh]
\caption{Summary of the outcomes of the LP model} \label{tab:LP}
\setlength{\tabcolsep}{3pt}
\begin{tabular}{|l|r|}
\hline
Number of Italian plants shut down& 1.704 \\
$\Delta$ Power (\%)& -65.46 \\
$\Delta$ Served population - Italy& +101,666,064 \\
$\Delta$ Served population - abroad& +3,522,406 \\
\hline
\end{tabular}
\end{table}

\noindent A brief comment concerning the comparison of the results reported in Tables \ref{tab:MILP} and \ref{tab:LP}.  The LP model has a significant impact on the Italian networks. Power reduction is greater than that produced by the MILP model ($65.46$\% vs. $53.23$\%) and, analogously, the number of shut-down Italian plants is greater (1,704 vs. 1,473). By contrast, the MILP model produces a greater number of incremental users of the Italian (133,824,096 vs. 101,666,064) and non-Italian networks (16,531,361 vs. 3,522,406). In both cases, we achieve the goals of reducing the power and increasing the service, but with an inverse ratio between power reduction and incremental users.

\medskip\noindent
In the most favorable case of the LP model, the annual energy used by the Italian FM networks would decrease, if the reduction obtained by the optimization were implemented, from $253$ to just over $152$ GWh ($40\%$ reduction) and the Italian FM "stove" would reduce its power from the current 28,800 KW to only 17,460 KW. Obviously, the saved power could be used elsewhere: Seen from the perspective of clean energy production, it is as if we had activated a new, small, "green" power plant of over $11$ MW.

\medskip\noindent
Finally, in Tables \ref{tab:it-networks} and \ref{abroad-networks} we provide some examples of the increment of the service areas both in Italy and in the bordering Administrations. In particular, we have listed the $20$ networks of Italian and non-Italian administrations with the greatest increase in users produced by the MILP and LP models. The Italian networks are anonymized, but is evident that the increment of users is uniformly distributed, which is another feature of our model.

\section{Simulation results}\label{sec:Simulation}

\medskip\noindent
Figures \ref{fig3} and \ref{fig4} show the service maps produced by the simulation tool. The color code is described in Section \ref{sec:Definitions} and the maps show the effect of implementing the solution produced by the LP model for a typical Italian network. In this case, Italy is the Affected Administration and all bordering Administrations are Interfering.

\Figure[hhh](topskip=0pt, botskip=0pt, midskip=0pt)[width=0.9\linewidth]{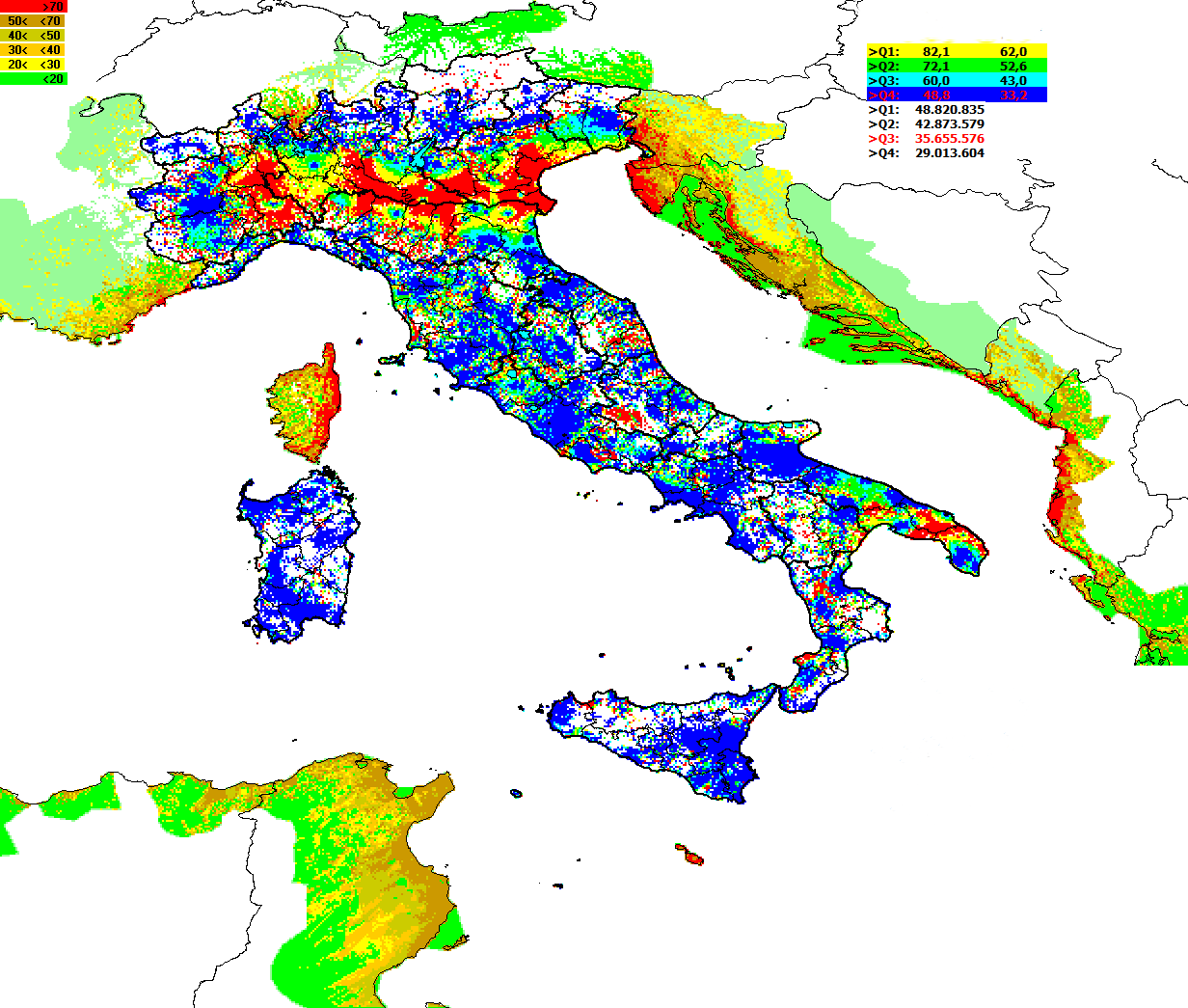} {Italian FM network:Current Scenario.\label{fig3}}

\Figure[hhh](topskip=0pt, botskip=0pt, midskip=0pt)[width=0.9\linewidth]{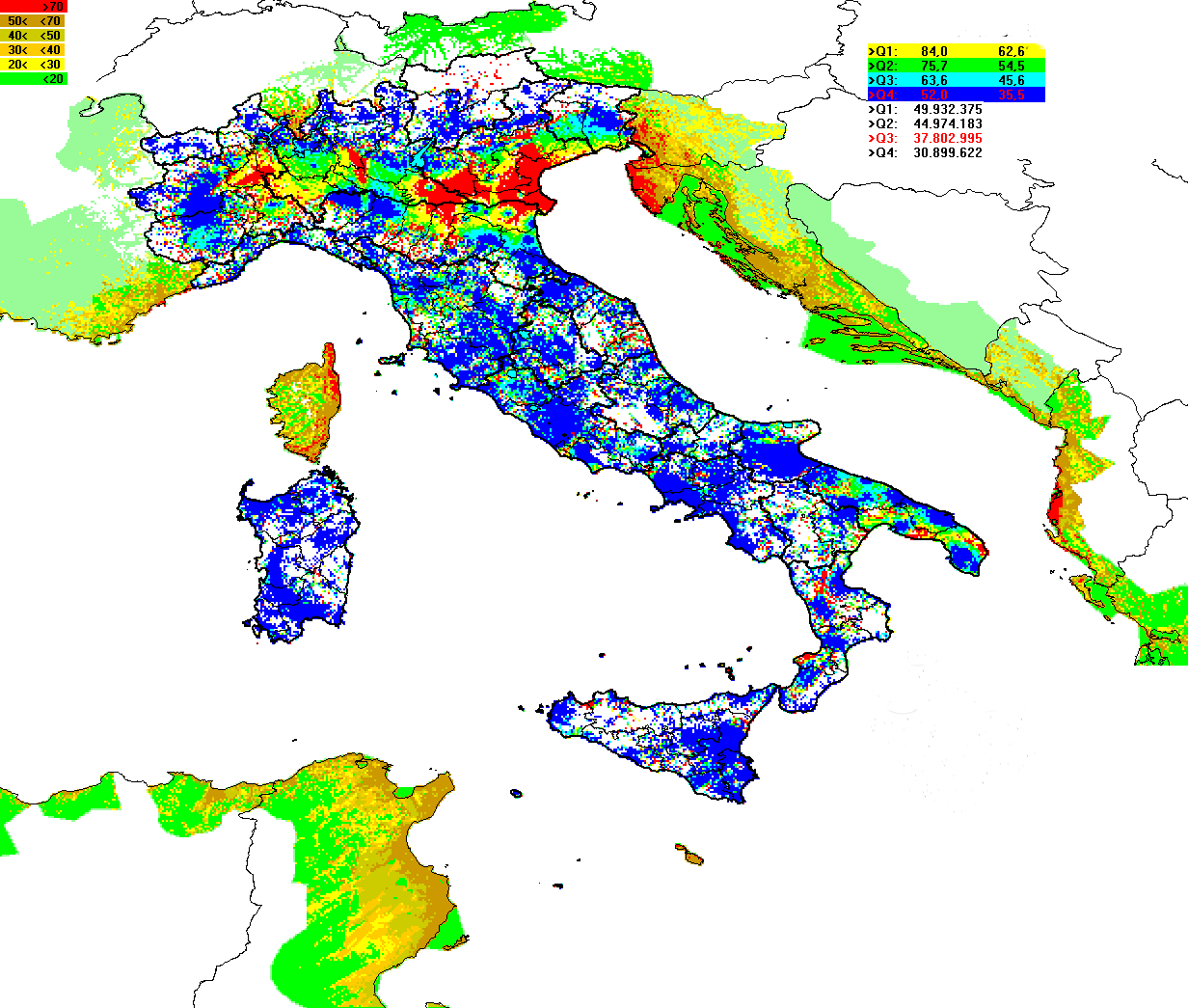} {Italian FM network: Post Power Reduction.\label{fig4}}

Fig.\ref{fig3} depicts the current situation of the Quality of Service for Italian receiving points and the interference produced by the Italian networks in the bordering Administrations, whereas Fig. \ref{fig4} depicts the resulting scenario after the application of the power reduction computed by the LP Model ($y$ variables). The maps are anonymized, but is clearly visible that the population served at quality Q3 ranges from $35,655,676$ to $37,802,995$ with more than $2$ million incremental users. Incremental users computed using the LP model and reported in Table \ref{tab:it-networks} are slightly greater than $1.6$ million (half a million more!). 

\begin{table}
\caption{Service improvement for some (anonymous) Italian networks}\label{tab:it-networks}
\begin{tabular}{|l|r|r|r|r|r|}
\hline
\#& Pop. Now& $\Delta$ Pop. MILP & $\Delta$ Pop. LP \\
\hline
$1$& 51,131,899& 1,151,046& 906,729 \\
$2$& 48,859,575& 1,586,859& 1,030,842 \\
$3$& 48,348,080& 1,998,243& 1,476,221 \\
$4$& 42,622,905& 2,028,532& 1,596,385 \\
$5$& 41,769,484& 3,433,057& 2,327,596 \\
$6$& 41,429,485& 2,703,513& 2,036,093 \\
$7$& 37,227,052& 1,334,123& 896,442 \\
$8$& 37,152,381& 2,975,092& 2,423,119 \\
$9$& 36,527,932& 4,283,262& 3,365,787 \\
$10$& 36,014,797& 2,475,593& 1,848,234 \\
$11$& 35,713,731& 2,361,952& 1,655,042 \\
$12$& 35,157,455& 2,259,645& 1,503,116 \\
$13$& 34,403,832& 1,779,844& 1,386,900 \\
$14$& 33,733,807& 2,079,745& 1,627,538 \\
$15$& 29,105,075& 2,375,934& 1,397,902 \\
$16$& 27,773,242& 2,677,454& 2,286,051 \\
$17$& 26,871,897& 2,407,246& 1,809,111 \\
$18$& 26,870,571& 2,722,907& 2,301,353 \\
$19$& 24,035,047& 1,807,826& 1,522,963 \\
$20$& 13,931,454& 643,932& 416,682 \\
\hline
\end{tabular}
\end{table}

\begin{table}
\caption{Service improvements for some bordering Administrations}\label{abroad-networks}
\begin{tabular}{|l|l|r|r|r|r|}
\hline
\# & Adm & Ch & Pop. Now & $\Delta$ Pop. LP & $\Delta$ Pop. MILP \\
\hline
$1$& F & 94,1  &	523,681	&	116,148	&	116,148	\\
$2$& TUN & 93,9  &	1,612,963	&	86,373	&	86,373	\\
$3$& ALB & 92,8 &	541,257	&	72,439	&	72,439	\\
$4$& ALB & 100 &	49,029	&	72,351	&	72,351	\\
$5$& F & 96,8 &	606,811	&	71,129	&	71,444	\\
$6$& TUN & 92,7 &	1,010,310	&	67,131	&	16,722	\\
$7$& HRV & 92,6  &	68,604	&	62,366	&	61,885	\\
$8$& F & 100,7 &	826,403	&	60,314	&	59,960	\\
$9$& F & 106 &	91,187	&	55,155	&	38,930	\\
$10$& HRV & 96,6 	&	124,434	&	49,872	&	49,707	\\
$11$& TUN & 89,5 	&	599,781	&	49,620	&	49,620	\\
$12$& F & 95 	&	459,459	&	46,429	&	46,429	\\
$13$& ALB & 90,7 	&	65,628	&	40,652	&	40,498	\\
$14$& HRV & 104,5 &	228,472	&	39,436	&	38,788	\\
$15$& SUI & 93,6  &	59,785	&	39,062	&	39,062	\\
$16$& SUI & 101 	&	2,563	&	37,000	&	37,000	\\
$17$& ALB & 101,6 	&	1,073,895	&	36,239	&	36,239	\\
$18$& F & 97,1 	&	188,315	&	35,987	&	35,987	\\
$19$& F& 101,2 &	187,857	&	35,230	&	35,230	\\
$20$& HRV & 99,3 &	446,505	&	34,643	&	35,112	\\
\hline
\end{tabular}
\end{table}

Hence we see that, even though a little obfuscated by the constraints of anonymity, there is an interesting effect. The simulated population covered by the optimized network is often greater than that computed by the optimization algorithm. This is due to the fact (among others) that the simulator does not fix the server $t_{ra}$ for a given network $a \in A$ and a given receiving point $r \in R$, but computes the SINR for all the potential servers of the network $a$ in $r$ and chooses the maximum. Hence, the power reduction preserves the original servers but also "promotes" server transmitters that were impaired by interference in the current scenario. 

\medskip\noindent
The effect of the optimization is particularly noticeable in several critical spots; for example, the service in Milan and Trieste in Italy is greatly improved and the interference caused in Corsica, Malta, and Albania is strongly reduced. Note also the striking effect of Power Reduction in the Italian region of \textit{Puglia} (\textit{South East Adriatic}) and Albania: The area served at quality Q4 in Puglia is significantly extended (in the main cities of Foggia, Bari and Brindisi) while the interference in Albania is strongly reduced. A "win-win" scenario.

\section{Conclusions}

\noindent In this paper we used a mixed-integer programming model to analyze the relation between the transmission power and population coverage in an FM radio infrastructure. The case study is based on an extensive dataset involving all Italian FM transmitters and all transmitters interfering with the Italian networks present in the ITU Master Register. Experimental results show that reducing the transmission power can lead to an even better coverage of the population, thanks to the reduction in interference. If implemented by regulatory bodies, these counterintuitive outcomes could bring significant benefits in terms of the quality of service, energy consumption, and electromagnetic pollution.

\section*{Acknowledgements}
\noindent The authors wish to thank Antonio Provenzano, Marco Ricchiuti, and Massimo Brienza of AGCOM (Italian Authority for Telecommunications) for useful and insightful discussions and for providing the publicly available data of their Register in a manageable form.

\EOD

\end{document}